\documentclass[
  10pt,
  a4paper,
]{article}

\usepackage[T1]{fontenc}
\usepackage[utf8]{inputenc}
\usepackage{textcomp}
\usepackage[english]{babel}

\usepackage{layaureo}

\usepackage{amsmath,amssymb,amsopn,amsthm}

\usepackage{indentfirst}

\usepackage{etoolbox}
\usepackage[
  mincrossrefs=99,
  style=numeric-comp,
  bibstyle=siam,
  sorting=nyt,
  backref=true,
  backend=biber,
  url=true,
  doi=true,
  maxbibnames=99,
  maxnames=99,
  giveninits=true,
  useprefix=true,
  date=year,
  autolang=other
]{biblatex}
\DeclareCiteCommand{\citeauthornsc}
  {%
   \boolfalse{citetracker}%
   \boolfalse{pagetracker}%
   \usebibmacro{prenote}}
  {\ifciteindex
     {\indexnames{labelname}}
     {}%
   \printnames{labelname}}
  {\multicitedelim}
  {\usebibmacro{postnote}}
\DeclareSourcemap{
  \maps[datatype=bibtex]{
    \map{
      \step[fieldsource=doi,final]
      \step[fieldset=url,null]
    }
  }
}

\usepackage[babel]{csquotes}
\usepackage{graphicx}
\usepackage{xcolor}

\definecolor{webgreen}{rgb}{0,.5,0}
\definecolor{webbrown}{rgb}{.6,0,0}
\definecolor{myblue}{rgb}{0,0.25,0.5}

\usepackage[
  pdfencoding=auto,
  colorlinks=true,
  linktocpage=true,
  pdfstartpage=1,
  pdfstartview=FitV,
  breaklinks=true,
  pdfpagemode=UseNone,
  pageanchor=true,
  pdfpagemode=UseOutlines,
  plainpages=false,
  bookmarksnumbered,
  bookmarksopen=true,
  bookmarksopenlevel=1,
  hypertexnames=true,
  pdfhighlight=/O,
  urlcolor=webbrown,
  linkcolor=myblue,
  citecolor=webgreen,
  pdfpagelabels,
]{hyperref}

\usepackage[osf,sc,noBBpl]{mathpazo} % don't use mathpazo's bb font (but amsmath's instead): solves PDF/A compliance problem
\usepackage{microtype}

\usepackage[capitalize,nameinlink]{cleveref}[0.19]
\crefname{section}{section}{sections}
\crefname{subsection}{subsection}{subsections}
\Crefname{figure}{Figure}{Figures}
\crefformat{equation}{\textup{#2(#1)#3}}
\crefrangeformat{equation}{\textup{#3(#1)#4--#5(#2)#6}}
\crefmultiformat{equation}{\textup{#2(#1)#3}}{ and \textup{#2(#1)#3}}
{, \textup{#2(#1)#3}}{, and \textup{#2(#1)#3}}
\crefrangemultiformat{equation}{\textup{#3(#1)#4--#5(#2)#6}}%
{ and \textup{#3(#1)#4--#5(#2)#6}}{, \textup{#3(#1)#4--#5(#2)#6}}{, and \textup{#3(#1)#4--#5(#2)#6}}
\Crefformat{equation}{#2Equation~\textup{(#1)}#3}
\Crefrangeformat{equation}{Equations~\textup{#3(#1)#4--#5(#2)#6}}
\Crefmultiformat{equation}{Equations~\textup{#2(#1)#3}}{ and \textup{#2(#1)#3}}
{, \textup{#2(#1)#3}}{, and \textup{#2(#1)#3}}
\Crefrangemultiformat{equation}{Equations~\textup{#3(#1)#4--#5(#2)#6}}%
{ and \textup{#3(#1)#4--#5(#2)#6}}{, \textup{#3(#1)#4--#5(#2)#6}}{, and \textup{#3(#1)#4--#5(#2)#6}}
\crefdefaultlabelformat{#2\textup{#1}#3}

\crefname{chapter}{chapter}{chapters}
\crefname{appendix}{appendix}{appendices}
\crefname{subappendix}{section}{sections}
\Crefname{subappendix}{Section}{Sections}
\crefname{page}{page}{pages}

\theoremstyle{plain}
\newtheorem{theorem}{Theorem}
\newtheorem{lemma}[theorem]{Lemma}
\newtheorem{corollary}[theorem]{Corollary}
\theoremstyle{definition}

\newtheorem*{conjecture}{Conjecture}
\theoremstyle{remark}
\newtheorem{remarkx}[theorem]{Remark}
\newenvironment{remark}
  {\pushQED{\qed}\remarkx}
  {\popQED\endremarkx}

\usepackage[font=small,labelfont=bf,textfont=it,margin=1em,labelsep=period]{caption}

\newcommand{\dd}{\mathop{}\!\mathrm{d}}

\newcommand{\ii}{\mathrm{i}}

\renewcommand{\restriction}{\raise-.5ex\hbox{\ensuremath{\upharpoonright}}}

\usepackage{mathtools}
\DeclarePairedDelimiter{\abs}{\lvert}{\rvert}

\title{Spectra of evolution operators of a class of neutral renewal equations: theoretical and numerical aspects}
\author{Dimitri Breda, Davide Liessi, Sjoerd M. Verduyn Lunel}
\date{19 May 2023}

\begin{document}

\maketitle

\begin{abstract}
In this work we begin a theoretical and numerical investigation on the spectra of evolution operators of neutral renewal equations, with the stability of equilibria and periodic orbits in mind.
We start from the simplest form of linear periodic equation with one discrete delay and fully characterize the spectrum of its monodromy operator.
We perform numerical experiments discretizing the evolution operators via pseudospectral collocation, confirming the theoretical results and giving perspectives on the generalization to systems and to multiple delays.
Although we do not attempt to perform a rigorous numerical analysis of the method, we give some considerations on a possible approach to the problem.


\smallskip
\noindent \textbf{Mathematics Subject Classification (2020):}
34K08, % Spectral theory of functional-differential operators
34K40, % Neutral functional-differential equations
37M99, % Approximation methods and numerical treatment of dynamical systems: None of the above but in this section
65Q10, % Numerical methods for difference equations
65Q20 % Numerical methods for functional equations
\end{abstract}

\section{Introduction}

Delays appear naturally in several phenomena, pertaining to, e.g., control theory, population dynamics and epidemics (see, e.g., \cite{Erneux2009,ArinoHbidAitDads2006}).
Indeed, many models are based on delay differential equations (DDEs) and renewal equations (REs).
Through their right-hand sides, DDEs prescribe the value at current time of the derivative of the unknown function, while REs prescribe the value at current time of the unknown function itself; in both cases these values depend on the unknown function at the current and past times.

The interest in delay equations has its roots in the 1930s and the theory of DDEs was developed starting from the 1940s; see the fundamental monographs \cite{DiekmannVanGilsVerduynLunelWalther1995,Hale1977,HaleVerduynLunel1993} and the references therein.
The sun-star framework of \cite{DiekmannVanGilsVerduynLunelWalther1995} was later partially extended to cover REs in \cite{DiekmannGettoGyllenberg2008,BredaLiessi2021}.

Numerical methods have been proposed to approximate the spectra of evolution operators for DDEs and REs or of the infinitesimal generator of their solution semigroups in the autonomous case (see \cite{BredaLiessi2018,BredaLiessi2020,BredaMasetVermiglio2005,BredaMasetVermiglio2012,BredaMasetVermiglio2015} and the references therein for pseudospectral collocation; see also \cite{BorgioliHajduInspergerStepanMichiels2020,LehotzkyInsperger2016,ButcherBobrenkov2011}).
They allow to study, e.g., the stability of equilibria and periodic orbits via the principle of linearized stability (see, e.g, \cite{DiekmannVanGilsVerduynLunelWalther1995} for DDEs and \cite{DiekmannGettoGyllenberg2008,BredaLiessi2021} for REs).

Recently, a new perturbation theory based on twin semigroups was proposed in \cite{DiekmannVerduynLunel2021}, where the authors derive the variation of constants formula not only for DDEs and REs but also for their neutral counterparts, characterized by unbounded perturbations of the trivial semigroup of evolution operators, as opposed to bounded ones.

Neutral delay equations have been used in mathematical models: as an example, neutral DDEs (NDDEs) often emerge from coupled oscillatory systems \cite{KyrychkoHogan2010} and neutral REs (NREs) from considering cohorts in cell populations \cite{DiekmannGettoNakata2016}.
In practice, NDDEs are typically characterized by the presence of delayed values of the derivative (of highest order) of the unknown function, while NREs typically involve discrete delay terms, as opposed to REs proper, which typically are integral equations.
Contrary to NDDEs \cite{Hale1977}, there are basically no theoretical results about NREs, except for those in \cite{DiekmannVerduynLunel2021} mentioned above (the interested reader may have a look also at \cite{DiekmannVanGils2000} for delay difference equations in discrete time).
Note that in light of \cite{DiekmannVerduynLunel2021} a theory similar to that of REs proper is expected to hold true, in particular with respect to the principle of linearized stability for equilibria and periodic orbits, as noted in \cite[section 13]{DiekmannVerduynLunel2021}.

To the best of our knowledge no specific numerical tools are available for the stability analysis of NREs, although DDE-BIFTOOL%
\footnote{\url{http://ddebiftool.sourceforge.net/}}
\cite{EngelborghsLuzyaninaRoose2002,
SieberEngelborghsLuzyaninaSamaeyRoose2014}
can analyze both NDDEs and delay differential algebraic equations \cite{KrauskopfSieber2023}, which potentially results in support for NREs as well when setting the left-hand side of a delay differential algebraic equation to zero.

The lack of theoretical and numerical tools hinders the adoption of NREs in mathematical models, while with the scarcity of models based on NREs the development of the former is at risk of missing a strong motivation.
This work aims at breaking this vicious cycle.
We are interested in particular in the stability of equilibria and periodic orbits of NREs and, although the relevant theory is currently lacking, we adopt the approach of studying the spectra of the evolution operators of the corresponding linearized equations.
Given the consolidated tradition of the cited numerical methods for delay equations, it seems natural to extend to NREs the pseudospectral collocation of \cite{BredaLiessi2018,BredaLiessiVermiglio2022}.
Our investigation develops both on the theoretical and on the numerical side, using the experiments as a guide in investigating and understanding the spectral theory of linear NREs, and using the proved theoretical results to validate the numerical approach.
The main theoretical result is the full characterization of the spectrum of the monodromy operator of a class of linear periodic NREs with one discrete delay and its decomposition in point, continuous and residual spectrum.
From the numerical point of view, instead, the proposed method works in general for any evolution operator of linear NREs, not necessarily periodic, with finite discrete and distributed delays.

The paper is organized as follows.
In \cref{sec:equation} we formulate the class of NREs for which we develop the theoretical results proved in \cref{sec:theory}.
Then, after briefly describing the discretization approach in \cref{sec:discretization}, we present in \cref{sec:experiments} a collection of numerical experiments exemplifying the theoretical results.
Finally, we give some perspectives on the generalization to systems of NREs with one delay in \cref{sec:system} and to scalar NREs with two delays in \cref{sec:2delays}.
Although we do not attempt here to perform a rigorous numerical analysis of the method, which is left to future work, some considerations are given in \cref{sec:conclusions}.

\section{Formulation of the problem}
\label{sec:equation}

To begin our investigation of evolution operators of NREs, with the perspective of studying their dynamic properties, we focus on the simplest linear NRE with constant or periodic coefficients.
We thus consider the scalar linear NRE (or difference equation)
\begin{equation}
\label{nre}
x(t) = f(t) x(t - \tau)
\end{equation}
with $\tau > 0$ and $f\colon \mathbb{R} \to \mathbb{R}$ a periodic function of bounded variation, continuous from the right.
For simplicity, in this work we restrict to the case where either $f$ is constant or its minimal period is equal to the delay $\tau$.
According to \cite{DiekmannVerduynLunel2021}, the natural state space for \cref{nre} is the space%
\footnote{We call the state space $Y$ for uniformity of notation with \cite{DiekmannVerduynLunel2021}.}
$Y \coloneqq NBV([-\tau, 0]; \mathbb{R})$ of real-valued functions of bounded variation on $[-\tau, 0]$ which are continuous from the right and have value $0$ at~$0$.
The space $Y$ is a Banach space with the total variation norm.
The initial value problem (IVP) associated to \cref{nre} is
\begin{equation}
\label{nre-ivp}
\left\{
\begin{aligned}
& x(t) = f(t) x(t - \tau), \quad t \geq s, \\
& x_{s} = \phi \in Y,
\end{aligned}
\right.
\end{equation}
where $x_{t}$ is the standard notation for the segment of $x$ at $t$ defined as $x_{t}(\theta) \coloneqq x(t+\theta)$ for $\theta \in [-\tau, 0]$.
Observe that the IVP \cref{nre-ivp} admits a unique solution for each $\phi \in Y$ (e.g., use the method of steps \cite{Bellman1961}, which consists in solving the equation on consecutive intervals of length $\tau$, so that at each step the past is fully known).
We can thus consider the associated monodromy operator, which is the evolution operator%
\footnote{Recall in general that an evolution operator $U(t,s)$ of a dynamical system maps the state of the system at time $s$ to the state at time $t\geq s$.}
$U \colon Y \to Y$ advancing the state of the system by one period along the solution, i.e.,
\begin{equation}
\label{U}
U \phi = U x_{0} \coloneqq x_{\tau}(\cdot; \phi),
\end{equation}
where $x(\cdot; \phi)$ is the solution of \cref{nre-ivp} with initial value $\phi$ at $t = 0$ ($x_0 = \phi$ and $x_{\tau}$ refer again to segments of $x$).
If $f$ is constant, we still define the operator $U$ as the evolution operator advancing the state of the system by a time $\tau$: this is customary with DDEs \cite{BredaMasetVermiglio2012} and REs \cite{BredaLiessi2018} when the stability of equilibria is investigated through the spectrum of evolution operators.

\section{Theoretical results}
\label{sec:theory}

In this \namecref{sec:theory} we fully characterize the spectrum of the resolvent set of \cref{nre} by assuming that the hypotheses on $f$ described in \cref{sec:equation} hold.

\begin{lemma}
\label{lem:Uf}
$U$ is the multiplication operator by $f$, i.e.,
$U\phi = f \phi$ for each $\phi \in Y$.
\end{lemma}
\begin{proof}
For $\theta \in [-\tau, 0]$ we have
\begin{equation*}
(U\phi)(\theta) = x_{\tau}(\theta;\phi) = x(\tau+\theta;\phi) = f(\tau+\theta) x(\theta;\phi) = f(\theta) \phi(\theta).
\qedhere
\end{equation*}
\end{proof}

In order to study the spectrum of $U$, we consider its complexification (recall indeed that the spectrum of $U$ is defined as the spectrum of its complexification).
With $Y_{\mathbb{C}} \cong NBV([-\tau, 0]; \mathbb{C})$, the complexified operator $U_{\mathbb{C}} \colon Y_{\mathbb{C}} \to Y_{\mathbb{C}}$ acts separately on the real and imaginary parts of $\phi = \phi_{\Re} + \ii \phi_{\Im} \in Y_{\mathbb{C}}$, i.e., $U_{\mathbb{C}}\phi \coloneqq U\phi_{\Re} + \ii U\phi_{\Im}$.
In the following, let $I_{Y}$ and $I_{Y_{\mathbb{C}}}$ be, respectively, the identity operators on $Y$ and $Y_{\mathbb{C}}$.

\begin{lemma}
\label{lem:Ulbij}
Let $\lambda = \alpha + \ii\beta \in \mathbb{C}$.
The operator $U_{\mathbb{C}} - \lambda I_{Y_{\mathbb{C}}}$ is bijective if and only if for each pair $(\psi_{\Re}, \psi_{\Im}) \in Y \times Y$ there exists a unique pair $(\phi_{\Re}, \phi_{\Im}) \in Y \times Y$ such that
\begin{equation}
\label{eq:Ulbijcond}
\left\{
\begin{aligned}
& (f - \alpha)\phi_{\Re} + \beta\phi_{\Im} = \psi_{\Re}, \\
& (f - \alpha)\phi_{\Im} - \beta\phi_{\Re} = \psi_{\Im}.
\end{aligned}
\right.
\end{equation}
The operator $U_{\mathbb{C}} - \lambda I_{Y_{\mathbb{C}}}$ is injective if and only if
\begin{equation*}
\left\{
\begin{aligned}
& (f - \alpha)\phi_{\Re} + \beta\phi_{\Im} = 0, \\
& (f - \alpha)\phi_{\Im} - \beta\phi_{\Re} = 0
\end{aligned}
\right.
\end{equation*}
for some $(\phi_{\Re}, \phi_{\Im}) \in Y \times Y$ implies that $\phi_{\Re} = \phi_{\Im} = 0$.
\end{lemma}
\begin{proof}
Observe that by separating the real and imaginary parts of $\phi, \psi \in Y_{\mathbb{C}}$, the equation $(U_{\mathbb{C}} - \lambda I_{Y_{\mathbb{C}}}) \phi = \psi$ is equivalent to $U\phi_{\Re} + \ii U\phi_{\Im} - \lambda\phi_{\Re} - \ii\lambda\phi_{\Im} = \psi_{\Re} + \ii\psi_{\Im}$, i.e., thanks to \cref{lem:Uf}, $f\phi_{\Re} + \ii f\phi_{\Im} - \lambda\phi_{\Re} - \ii\lambda\phi_{\Im} = \psi_{\Re} + \ii\psi_{\Im}$.
As for the second part, recall that a linear operator is injective if and only if its kernel is trivial.
\end{proof}

\begin{theorem}
\label{th:sU}
The spectrum and resolvent set of $U_{\mathbb{C}}$ are, respectively,
\begin{equation*}
\sigma(U_{\mathbb{C}}) = \overline{f(\mathbb{R})},
\qquad
\rho(U_{\mathbb{C}}) = \mathbb{C}\setminus\overline{f(\mathbb{R})}.
\end{equation*}
\end{theorem}
\begin{proof}
We proceed by subsequently proving that $f(\mathbb{R}) \subset \sigma(U_{\mathbb{C}})$, that $\overline{f(\mathbb{R})} \setminus f(\mathbb{R}) \subset \sigma(U_{\mathbb{C}})$, that $\mathbb{R} \setminus \overline{f(\mathbb{R})} \subset \rho(U_{\mathbb{C}})$ and that $\mathbb{C} \setminus \mathbb{R} \subset \rho(U_{\mathbb{C}})$.

Let $\lambda \in f(\mathbb{R}) \subset \mathbb{R}$ and let $\hat{\theta} \in [-\tau,0)$ be such that $f(\hat{\theta}) = \lambda$.
Then, recalling \cref{lem:Uf,lem:Ulbij}, for all $\phi \in Y_{\mathbb{C}}$ we have
\begin{equation}
\label{eq:Ulzero}
((U_{\mathbb{C}} - \lambda I_{Y_{\mathbb{C}}})\phi)(\hat{\theta}) = (f(\hat{\theta}) - \lambda)\phi(\hat{\theta}) = 0,
\end{equation}
which implies that $U_{\mathbb{C}} - \lambda I_{Y_{\mathbb{C}}}$ is not surjective (any $\psi \in Y_{\mathbb{C}}$ such that $\psi(\hat{\theta}) \neq 0$ has no inverse image).
Hence $f(\mathbb{R}) \subset \sigma(U_{\mathbb{C}})$.

Let $\lambda \in \overline{f(\mathbb{R})} \setminus f(\mathbb{R}) \subset \mathbb{R}$.
Then there exists a sequence $\{\lambda_{n}\}_{n\in\mathbb{N}}$ in $f(\mathbb{R})$ such that $\lambda_n \to \lambda$, and a sequence $\{\theta_{n}\}_{n\in\mathbb{N}} \subset [-\tau, 0)$ (thanks to the periodicity of $f$) such that $f(\theta_{n}) = \lambda_{n}$.
The sequence $\{\theta_{n}\}_{n\in\mathbb{N}}$ is bounded, hence by the Bolzano--Weierstrass theorem it has a subsequence $\{\theta_{n_{m}}\}_{m\in\mathbb{N}}$ such that $\theta_{n_{m}} \to \hat{\theta}$ for some $\hat{\theta} \in [-\tau, 0]$.
Observe that thanks to \cref{lem:Uf,lem:Ulbij} for every $\phi \in Y$
\begin{equation}
\label{eq:fltozero}
\begin{aligned}
((U-\lambda I_{Y})\phi)(\theta_{n_{m}}) &= (f(\theta_{n_{m}}) - \lambda)\phi(\theta_{n_{m}}) \\
&= (f(\theta_{n_{m}}) - \lambda_{n_{m}} + \lambda_{n_{m}} - \lambda)\phi(\theta_{n_{m}}) \\
&= (\lambda_{n_{m}} - \lambda)\phi(\theta_{n_{m}}) \to 0,
\end{aligned}
\end{equation}
since $\lambda_{n_{m}} \to \lambda$ and $\phi$ has bounded variation and is thus bounded.
Let $\psi \in Y$ be continuous and not null at $\hat{\theta}$.
Thus there exists $\delta > 0$ such that for every $\theta$ with $0 < \abs{\theta-\hat{\theta}} < \delta$ we have $\abs{\psi(\theta)-\psi(\hat{\theta})} < \abs{\psi(\hat{\theta})}/2$, i.e., $\abs{\psi(\theta)} > \abs{\psi(\hat{\theta})}/2$.
From \eqref{eq:fltozero} for each $\phi \in Y$ there exists $M \in \mathbb{R}$ such that $\abs{((U-\lambda I_{Y})\phi)(\theta_{n_{m}})} <  \abs{\psi(\hat{\theta})}/2$ for every $m > M$.
Choosing $m > M$ such that $0 < \abs{\theta_{n_{m}} - \hat{\theta}} < \delta$ (observe that since $\lambda \not\in f(\mathbb{R})$ the sequence $\theta_{n_{m}}$ can be equal to $\hat{\theta}$ only for a finite number of values of $m$), we obtain that $\abs{((U-\lambda I_{Y})\phi)(\theta_{n_{m}})} <  \abs{\psi(\hat{\theta})}/2 < \abs{\psi(\theta_{n_{m}})}$ and $((U-\lambda I_{Y})\phi)(\theta_{n_{m}}) \neq \psi(\theta_{n_{m}})$.
Hence neither $U - \lambda I_{Y}$ nor $U_{\mathbb{C}} - \lambda I_{Y_{\mathbb{C}}}$ (consider $\psi+i0 \in Y_{\mathbb{C}}$) are surjective.
This implies that $\overline{f(\mathbb{R})} \setminus f(\mathbb{R}) \subset \sigma(U_{\mathbb{C}})$.

Let $\lambda \in \mathbb{R} \setminus \overline{f(\mathbb{R})}$.
There exists a neighborhood of $\lambda$ contained in $\mathbb{R} \setminus \overline{f(\mathbb{R})}$ and $f-\lambda$ is bounded away from $0$.
Let $\psi \in Y$.
Since both $\psi$ and $f-\lambda$ have bounded variation and are continuous from the right, and since $\psi(0)=0$, the function $\phi \coloneqq \psi / (f-\lambda)$ is in $Y$, thus $U - \lambda I_{Y}$ is bijective, and thanks to \cref{lem:Ulbij} also $U_{\mathbb{C}} - \lambda I_{Y_{\mathbb{C}}}$ is.
Hence $\mathbb{R} \setminus \overline{f(\mathbb{R})} \subset \rho(U_{\mathbb{C}})$.

Let $\lambda = \alpha + \ii\beta \in \mathbb{C} \setminus \mathbb{R}$, i.e., with $\beta \neq 0$.
From \cref{lem:Ulbij} the operator $U_{\mathbb{C}} - \lambda I_{Y_{\mathbb{C}}}$ is bijective if and only if for each pair $(\psi_{\Re}, \psi_{\Im}) \in Y \times Y$ there exists a unique pair $(\phi_{\Re}, \phi_{\Im}) \in Y \times Y$ such that \cref{eq:Ulbijcond} holds.
Solving for $\phi_{\Re}$ and $\phi_{\Im}$ we obtain
\begin{equation*}
\left\{
\begin{aligned}
& \phi_{\Re} = \frac{(f - \alpha)\psi_{\Re} - \beta\psi_{\Im}}{\beta^{2}+(f-\alpha)^{2}}, \\
& \phi_{\Im} = \frac{(f - \alpha)\psi_{\Im} + \beta\psi_{\Re}}{\beta^{2}+(f-\alpha)^{2}}.
\end{aligned}
\right.
\end{equation*}
Observe that $\beta^{2}+(f-\alpha)^{2}$ is bounded away from $0$, has bounded variation and is continuous from the right, so $\phi_{\Re}, \phi_{\Im} \in Y$, proving the bijectivity of $U_{\mathbb{C}} - \lambda I_{Y_{\mathbb{C}}}$.
Thus $\mathbb{C} \setminus \mathbb{R} \subset \rho(U_{\mathbb{C}})$.
\end{proof}

\begin{theorem}
The point spectrum of $U_{\mathbb{C}}$ is
\begin{equation*}
\begin{aligned}
\sigma_{p}(U_{\mathbb{C}})
&
= \{\lambda \in \sigma(U_{\mathbb{C}}) \mid (f^{-1}(\lambda))^{\circ} \neq \emptyset\}
\\
&
= \{\lambda \in \sigma(U_{\mathbb{C}}) \mid f^{-1}(\lambda) \text{ contains a segment}\}.
\end{aligned}
\end{equation*}
\end{theorem}
\begin{proof}
Let $\lambda \in \sigma(U_{\mathbb{C}}) \subset \mathbb{R}$.
Thanks to \cref{lem:Ulbij}, $U_{\mathbb{C}} - \lambda I_{Y_{\mathbb{C}}}$ is not injective if and only if there exist $\phi_{\Re}, \phi_{\Im} \in Y$, not both null, such that $(f-\lambda)\phi_{\Re} = 0$ and $(f-\lambda)\phi_{\Im} = 0$.
It is enough to consider the existence of one $\phi \in Y \setminus \{0\}$ such that $(f-\lambda)\phi = 0$.
This equality holds if and only if $f(\theta) = \lambda$ or $\phi(\theta) = 0$ for each $\theta \in [-\tau, 0]$.
Since $\phi \neq 0$, there exists $\hat{\theta} \in [-\tau, 0)$ such that $\phi(\hat{\theta}) \neq 0$ (recall the normalization condition $\phi(0) = 0$).
Since $\phi$ is continuous from the right, there is a right neighborhood of $\hat{\theta}$ on which $\phi$ is not null: on this neighborhood $f$ must assume the constant value $\lambda$.
Hence $U_{\mathbb{C}} - \lambda I_{Y_{\mathbb{C}}}$ is not injective if and only if there exists an interval in $[-\tau, 0]$ such that $f$ is constantly $\lambda$ on that interval, i.e., $f^{-1}(\lambda)$ contains a segment, or, in other words, the interior of $f^{-1}(\lambda)$ is not empty.
\end{proof}

\begin{corollary}
$\sigma_{p}(U_{\mathbb{C}}) \subset f(\mathbb{R})$.
\end{corollary}

\begin{theorem}
The continuous spectrum $\sigma_{c}(U_{\mathbb{C}})$ of $U_{\mathbb{C}}$ is empty.
\end{theorem}
\begin{proof}
Recalling \cref{eq:Ulzero} for some $\lambda$ and $\hat{\theta}$ with $\lambda = f(\hat{\theta}) \in \mathbb{R}$, consider $\psi \in Y_{\mathbb{C}}$ such that $\psi(\hat{\theta}) \neq 0$.
For each $\phi \in Y_{\mathbb{C}}$, let $\psi_{\phi} \coloneqq \psi - (U_{\mathbb{C}} - \lambda I_{Y_{\mathbb{C}}})\phi \in Y_{\mathbb{C}}$.
We have $\psi_{\phi}(0) = 0$ and, recalling \cref{lem:Uf}, $\psi_{\phi}(\hat{\theta}) = \psi(\hat{\theta})$, so the total variation of $\psi_{\phi}$ is at least $\abs{\psi(\hat{\theta})}$.
Hence no sequence of functions in $Y_{\mathbb{C}}$ (and thus in the range of $U_{\mathbb{C}} - \lambda I_{Y_{\mathbb{C}}}$) can converge in total variation norm to $\psi$, i.e., $f(\mathbb{R}) \cap \sigma_{c}(U_{\mathbb{C}}) = \emptyset$.

Consider now $\lambda \in \overline{f(\mathbb{R})} \setminus f(\mathbb{R}) \subset \mathbb{R}$.
As in the proof of \cref{th:sU}, we can construct a sequence $\{\theta_{n_{m}}\}_{m\in\mathbb{N}}$ in $[-\tau, 0)$ with some limit $\hat{\theta} \in [-\tau, 0]$ such that $f(\theta_{n_{m}})\to \lambda$ and \cref{eq:fltozero} holds for each $\phi \in Y$.
We can choose $\psi \in Y$ continuous and not null at $\hat{\theta}$ and obtain $\delta > 0$ such that for every $\theta$ with $0 < \abs{\theta-\hat{\theta}} < \delta$ we have $\abs{\psi(\theta)-\psi(\hat{\theta})} < \abs{\psi(\hat{\theta})}/4$, i.e., $\abs{\psi(\theta)} > 3\abs{\psi(\hat{\theta})}/4$.
From \eqref{eq:fltozero} for each $\phi \in Y$ there exists $M \in \mathbb{R}$ such that $\abs{((U-\lambda I_{Y})\phi)(\theta_{n_{m}})} <  \abs{\psi(\hat{\theta})}/4$ for every $m > M$.
We now choose $m > M$ such that $0 < \abs{\theta_{n_{m}} - \hat{\theta}} < \delta$ and obtain that $\abs{((U-\lambda I_{Y})\phi)(\theta_{n_{m}})} < \abs{\psi(\hat{\theta})}/4 < 3\abs{\psi(\hat{\theta})}/4 < \abs{\psi(\theta_{n_{m}})}$.
The total variation of $\psi - (U - \lambda I_{Y})\phi$ is thus at least $\abs{\psi(\hat{\theta})}/2$ and no sequence of functions in $Y$ can converge in total variation norm to $\psi$.
The same holds for the complexification and thus for the range of $U_{\mathbb{C}} - \lambda I_{Y_{\mathbb{C}}}$.
Then also $(\overline{f(\mathbb{R})} \setminus f(\mathbb{R})) \cap \sigma_{c}(U_{\mathbb{C}}) = \emptyset$.
\end{proof}

\begin{corollary}
The residual spectrum of $U_{\mathbb{C}}$ is
\begin{equation*}
\sigma_{r}(U_{\mathbb{C}}) = \overline{f(\mathbb{R})} \setminus \{\lambda \in \mathbb{C} \mid (f^{-1}(\lambda))^{\circ} \neq \emptyset\}.
\end{equation*}
\end{corollary}
\begin{proof}
Recall that $\sigma_{r}(U_{\mathbb{C}}) = \sigma(U_{\mathbb{C}}) \setminus (\sigma_{p}(U_{\mathbb{C}}) \cup \sigma_{c}(U_{\mathbb{C}}))$.
\end{proof}

\begin{theorem}
If $f \not \equiv 0$, the operators $U$ and $U_{\mathbb{C}}$ are not compact.
\end{theorem}
\begin{proof}
If there exists $[a, b) \subset [-\tau, 0]$ such that $f\restriction_{[a, b)}$ is constant and not null, define the spaces $Y_{1} \coloneqq NBV([-\tau, a]; \mathbb{C})$, $Y_{2} \coloneqq NBV([a, b]; \mathbb{C})$ and $Y_{3} \coloneqq NBV([b, 0]; \mathbb{C})$.
Then, considering the immersions given by prolongation with $0$, $Y \cong Y_{1} \oplus Y_{2} \oplus Y_{3}$.
Let $U_{i} \coloneqq U_{\mathbb{C}}\restriction_{Y_{i}}$ for $i \in \{1,2,3\}$ and observe that $U_{2}$ is a nonzero multiple of the identity, so it is not compact.
Since restrictions of compact operators are compact, $U_{\mathbb{C}}$ is not compact.
The same holds for $U$.

If $f$ has no constant pieces, thanks to the continuity from the right, the image of $f$ contains a segment and is thus uncountable.
Since a compact operator has countable spectrum and $\sigma(U) = \sigma(U_{\mathbb{C}})$, $U$ and $U_{\mathbb{C}}$ are not compact.
The same argument applies to the last case of $f$ having some constant pieces, but all of them being null, recalling that $f\not\equiv0$, so there exists $t\in\mathbb{R}$ such that $f(t) \neq 0$.
\end{proof}

\begin{remark}
For a study of multiplication operators and their spectra in the slightly different context of bounded variation functions (without normalization and continuity from the right), see \cite{AstudilloVillalbaRamosFernandez2017}.
In particular, similarly to \cref{th:sU}, it is shown in \cite[Theorem 13]{AstudilloVillalbaRamosFernandez2017} that the spectrum of the multiplication operator by $f$ is the closure of the range of $f$.
However, when the condition of continuity from the right is imposed on the function space, some results therein do not hold up (e.g., \cite[Proposition 2]{AstudilloVillalbaRamosFernandez2017}); moreover, most proofs need to be adapted or do not work anymore, more specifically those based on the construction of functions with specific values at a finite number of points and another value elsewhere, which are obviously not continuous from the right.
\end{remark}

\section{Discretization of the evolution operators}
\label{sec:discretization}

In \cref{sec:theory} it has been possible to obtain sharp theoretical results on the spectrum of evolution operators of NREs by restricting the analysis to the specific class \cref{nre} of scalar periodic NREs with a single discrete and constant delay and period equal to this delay.
From the numerical point of view, instead, the pseudospectral collocation proposed in \cite{BredaLiessi2018} for REs proper (and then extended in several directions in \cite{BredaLiessi2020,BredaLiessiVermiglio2022,BredaLiessiVermiglio}) is rather general.
Indeed, at least from the implementation point of view, the method of \cite{BredaLiessi2018} can be applied to NREs as well.
We thus recall its essentials in the following by referring to the evolution family $\{U(t,s)\}_{t\geq s}$ associated to the IVP
\begin{equation}
\label{nre-gen-ivp}
\left\{
\begin{aligned}
& x(t)=F(t,x_{t}), \quad t\geq s, \\
& x_{s}=\phi \in Y,
\end{aligned}
\right.
\end{equation}
assuming the latter to be well-posed and that $F \colon \mathbb{R} \times Y \to \mathbb{R}$ is linear in the second argument.%
\footnote{For simplicity we restrict to a scalar equation; extending to any (finite) dimension is straightforward.}
In particular, and without loss of generality, we describe the discretization of $U\coloneqq U(s+h,s)$ given any $s\in\mathbb{R}$ and any $h>0$.

The first step is a reformulation of the monodromy operator, which is convenient for discretizing the operator, approximating its spectrum and (at least in the case of DDEs \cite{BredaMasetVermiglio2012} and REs \cite{BredaLiessi2018}) proving the convergence of the method.
We define the auxiliary function spaces $Y^{+} \coloneqq NBV([0, h]; \mathbb{R})$ and $Y^{\pm} \coloneqq NBV([-\tau, h]; \mathbb{R})$, the operator $V \colon Y \times Y^{+} \to Y^{\pm}$ as
\begin{equation*}
V(\phi, w) (t) \coloneqq
\begin{cases}
w(t), & t \in (0, h], \\
\phi(t), & t \in [-\tau, 0],
\end{cases}
\end{equation*}
and the operator $\mathcal{F}_{s} \colon Y^{\pm} \to Y^{+}$ as
\begin{equation*}
\mathcal{F}_{s}u(t) \coloneqq F(s+t, u_t), \quad t \in [0, h].
\end{equation*}
The operator $U \colon Y \to Y$ defined in \cref{U} is then reformulated as
\begin{equation}
\label{U-reform1}
U\phi = V(\phi, w^{\ast})_{h},
\end{equation}
where $w^{\ast}$ is the unique solution of the fixed point equation
\begin{equation}
\label{U-reform2}
w = \mathcal{F}_s V(\phi, w).
\end{equation}
Indeed, it is clear that solutions of \cref{U-reform2} correspond to solutions on $[0, h]$ of \cref{nre-gen-ivp}, which we assumed to be well-posed.

Observe that the operator $\mathcal{F}_s$ consists in the application of the right-hand side of the equation (with a time-shift of $s$ to keep working on the intervals $[-\tau,0]$ and $[0, h]$), while the operator $V$ represents the type of equation and describes how the solution is constructed from the initial value and the output of the right-hand side.
Thus, for REs proper the same construction as above is applied to spaces of $L^1$ functions, while for DDEs, besides using spaces of continuous functions, the definition of $V$ is different, having $w$ the role of the derivative of the solution.

The pseudospectral collocation technique of \cite{BredaMasetVermiglio2012,BredaLiessi2018} consists in applying the monodromy operator to polynomials interpolating the functions in the state space and in the auxiliary spaces.
Functions of $Y$ and $Y^{+}$ are thus represented by the vectors of their values at fixed sets of nodes in $[-\tau, 0]$ and $[0, h]$ (e.g., Chebyshev nodes).
The dimension of these vectors depends on the degree $M$ of the polynomials and, in case of a piecewise approach, on the number $L$ of pieces.
Let $Y^{+}_{L,M} \coloneqq \mathbb{R}^{LM+1}$ be the real vector space corresponding to $Y^+$, $R^{+} \colon Y^{+} \to Y^{+}_{L,M}$ be the restriction operator associating to a function the corresponding vector and $P^{+} \colon Y^{+}_{L,M} \to Y^{+}$ be the prolongation operator associating to a vector the polynomial interpolating its values.
More precisely, let $0 = c_0 < c_1 < \dots < c_{M+1} = 1$ form the chosen family of (normalized) collocation abscissae, let $0 = t_0 < t_1 < \dots < t_L = h$ be the endpoints of the pieces in $[0, h]$ and define the collocation nodes as $t_{i,j} \coloneqq t_i+c_j(t_{i+1}-t_i)$ for $i \in \{0, \dots, L-1\}$ and $j \in \{0, M+1\}$ and $t_{L,0}\coloneqq h$ (observe that $t_{i,M+1} = t_{i+1,0}$); then for $w \in Y^{+}$ the restriction $R^{+}w$ is the vector
\begin{equation}
\label{Rw}
\begin{aligned}
(w(t_{0,0}), \dots, w(t_{0,M}), {}&{ }w(t_{1,0}), \dots, w(t_{1,M}), \dots, \\
& w(t_{L-1,0}), \dots, w(t_{L-1,M}), w(t_{L,0})),
\end{aligned}
\end{equation}
while for $W \in Y^{+}_{L,M}$ (with indices in the same order as in \cref{Rw}) the prolongation $P^{+}W$ is the piecewise polynomial interpolating the values $W_{i,0}, \dots, W_{i,M}, W_{i+1,0}$ in $[t_i,t_{i+1}]$ for each $i \in \{0, \dots, L-1\}$.
The real vector space $Y_{L,M}$ corresponding to $Y$ and the corresponding restriction and prolongation operators $R$ and $P$, respectively, are defined similarly, with the pieces in $[-\tau,0]$ defined by shifting the pieces in $[0, h]$ by (multiples of) $-h$ and using the same collocation parameters $c_0, \dots, c_{M+1}$; some attention is needed near $-\tau$, see \cite{BredaLiessiVermiglio2022} for more details.

Given the reformulation \cref{U-reform1,U-reform2}, we discretize the operator $U$ as the finite-dimensional operator $U_{L,M} \colon Y_{L,M} \to Y_{L,M}$ defined as
\begin{equation*}
U_{L,M} \Phi \coloneqq R V(P \Phi, P^{+} W^{\ast})_{h},
\end{equation*}
where $W^{\ast} \in Y^{+}_{L,M}$ is a solution of the fixed point equation
\begin{equation}
\label{fixed-point-UMN}
W = R^{+} \mathcal{F}_{s} V(P \Phi, P^{+} W)
\end{equation}
for the given $\Phi \in Y_{L,M}$.
The operator $U_{L,M}$ can now be represented as a matrix (see \cite[appendix A]{BredaLiessi2018} for details), whose eigenvalues can be computed with standard methods (e.g., MATLAB's \texttt{eig} function).

We consider these eigenvalues as the approximations of the elements of the spectrum of $U_{\mathbb{C}}$ for $L,M \to +\infty$.
The precise meaning of the limit and of the convergence in this context needs some attention, as detailed below comparing REs proper and NREs.

In the case of REs proper, given a monodromy operator $U$ of a periodic equation, with reasonable regularity conditions on $F$ each nonzero element of $\sigma(U_{\mathbb{C}})$ is isolated \cite{BredaLiessi2021}.
It is proved in \cite{BredaLiessi2018} under suitable assumptions that for each element of $\sigma(U_{\mathbb{C}})$ as $M$ increases there are eigenvalues in $\sigma(U_{1,M})$ converging to it.
The order of convergence depends on the regularity of the eigenfunctions and it is infinite in $M$ if they are smooth: this is a typical phenomenon for pseudospectral methods, often called \emph{spectral accuracy} (see \cite[chapter 4]{Trefethen2000} and \cite[chapter 2]{Boyd2001}).
Moreover, the piecewise method ($L \geq 1$) exhibits a finite order of convergence in $L$, see \cite[section 4.1]{BredaLiessi2020} for details.

For NREs of the type \cref{nre}, the spectrum of a monodromy operator $U$, seen as a subset of $\mathbb{R}$, may have a nonempty interior: consider for instance $f(t) = \sin(2 \pi t / \tau)$, for which $\sigma(U_{\mathbb{C}}) = [-1, 1]$ according to \cref{th:sU}.
For the points on the boundary of $\sigma(U_{\mathbb{C}})\subset\mathbb{R}$ it may still make sense to consider the convergence in the familiar way.
For points in the interior, instead, defining a precise notion of convergence is a delicate issue and in this work we prefer to rely on the intuition%
\footnote{One possible idea is that in the limit (in some sense to be defined) the set $\sigma(U_{L,M})$ should become dense in $\sigma(U_{\mathbb{C}})$.}
of the reader, informed by some experimental observations in \cref{sec:experiments}.
Apart from this complication, the convergence analysis of \cite{BredaLiessi2018} relies on the Banach Perturbation Lemma and requires that the equation have a regularizing effect on the solution, which is lacking in the case of NREs: an alternative approach is thus necessary and will be the subject of future research, although some relevant considerations are given in \cref{sec:conclusions}.

In order for $U_{L,M}$ to be well-defined, the solution of \cref{fixed-point-UMN} needs to exist and be unique for every $\Phi \in Y_{L,M}$.
For REs proper, in \cite[section 4.2]{BredaLiessi2018} conditions are given to ensure this in the general case.
For NREs, one needs to either assume it or prove it.
In the specific case of \cref{nre}, thanks to the periodicity of $f$, the existence and uniqueness can be easily proved: indeed, we can observe that for each $i$ and each $j$ we have $W^{\ast}_{i,j} = (\mathcal{F}_{s}V(P\Phi,P^{+}W^{\ast}))(t_{i,j}) = f(t_{i,j}) V(P\Phi,P^{+}W^{\ast})(t_{i,j}-\tau) = f(t_{i,j})(P\Phi)(t_{i,j}-\tau)$, so given $\Phi$ the vector $W^{\ast}$ is uniquely determined.
The proof in \cite{BredaLiessi2018} is again based on the Banach Perturbation Lemma and on the regularizing effect of the equation: investigating alternative approaches for NREs left to future work also in this case.

As a final note, observe that the numerical method does not explicitly compute the solution, which is needed only in the theoretical formulation of the method: instead, it directly approximates the operator.

\section{Experiments}
\label{sec:experiments}

Computing the spectrum of the monodromy operator $U$ \cref{U} with different choices for $f$ using the method described in \cref{sec:discretization} confirms the validity of the theoretical results collected in \cref{sec:theory}.
In particular the computed spectrum approximates the set $\overline{f(\mathbb{R})} \cup \{0\}$, so it seems that the method approximates the whole spectrum and not only the point spectrum.
In this regard, we may consider the numerical method reliable as far as this work is concerned, even in the absence of a proof of convergence.

\begin{figure}
\centering
\includegraphics{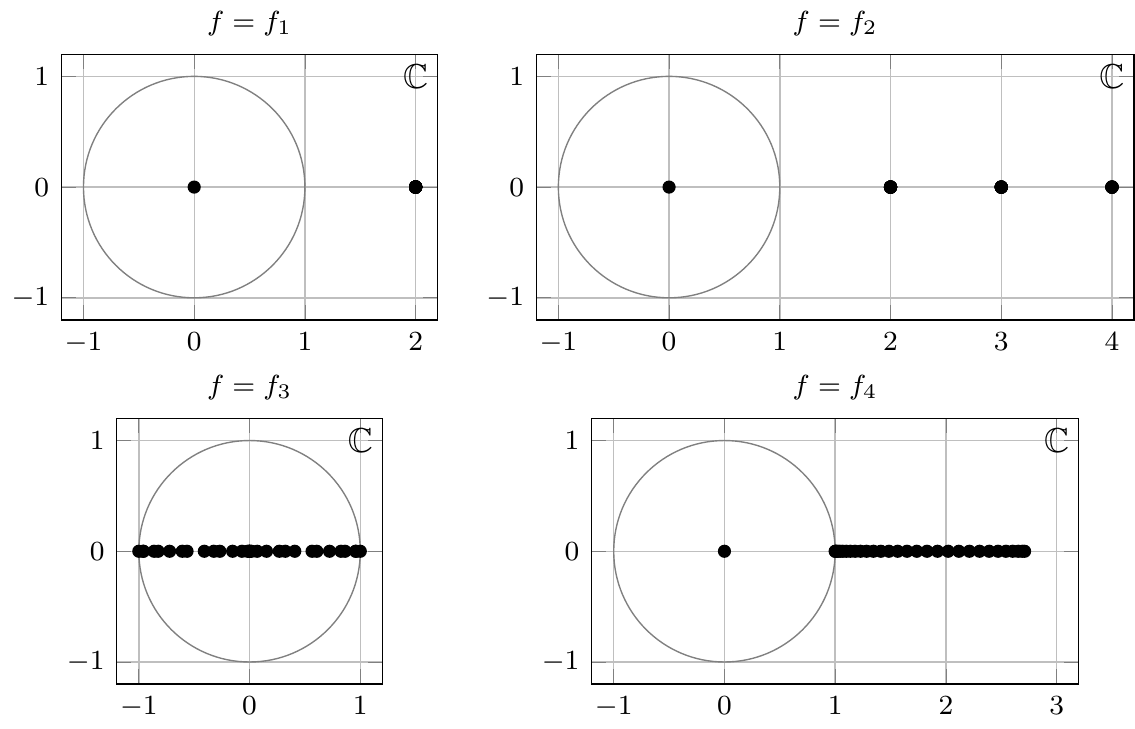}

\caption{Spectra of $U_{\mathbb{C}}$ with $f = f_i$ for $i \in \{1,2,3,4\}$ (see \cref{f1,f2,f3,f4}), computed with $L=1$ and $M=30$.}
\label{fig:spectra1}
\end{figure}

In the experiments we choose $\tau = 1$ without loss of generality.
\Cref{fig:spectra1} shows the computed spectra for the following choices of $f$:
\begin{alignat}{2}
\label{f1}
f_1(t) &\coloneqq 2, & \qquad & (\sigma(U_{\mathbb{C}}) = \{2\}) \\
\label{f2}
f_2(t) &\coloneqq
\begin{cases}
2 & \text{if $t-\lfloor t \rfloor \in [0, 0.3)$,} \\
3 & \text{if $t-\lfloor t \rfloor \in [0.3, 0.8)$,} \\
4 & \text{if $t-\lfloor t \rfloor \in [0.8, 1)$,}
\end{cases}
& \qquad & (\sigma(U_{\mathbb{C}}) = \{2, 3, 4\}) \\
\label{f3}
f_3(t) &\coloneqq \sin(2\pi t), & \qquad & (\sigma(U_{\mathbb{C}}) = [-1,1]) \\
\label{f4}
f_4(t) &\coloneqq \exp(t-\lfloor t \rfloor). & \qquad & (\sigma(U_{\mathbb{C}}) = [1, e])
\end{alignat}
The spectra are computed with $L=1$ and $M = 30$.

Unfortunately, in general a numerical method based on a matrix approximation has no means to discriminate between the point, continuous and residual spectra of an operator; therefore we can only experimentally verify \cref{th:sU} and not the rest of the results in \cref{sec:theory}.
We also cannot really observe the difference between a set and its closure numerically.
However, the results obtanied with $f = f_4$ suggest that the approximated spectrum actually is the closure of the image of $f$: indeed, $f_4(\mathbb{R}) = [1,e)$, but one of the approximated eigenvalues appears to converge to $e$ (with order $2$).

As observed above, the computed spectrum always contains $0$, even when it is not in the spectrum of $U_{\mathbb{C}}$.
This can be explained by observing that the operator on the state space $Y$ corresponding to the finite-dimensional operator $U_{L,M}$ has finite rank and is thus compact; the method actually computes the spectrum of that finite-rank compact operator, which has $0$ in the spectrum.

As we noted in \cref{sec:discretization}, unless we further specify the meaning of the convergence, it makes sense to talk about the order of convergence only for isolated points in the spectrum and, possibly, for points on its boundary (considering it as a subset of $\mathbb{R}$).

\begin{figure}
\centering
\includegraphics{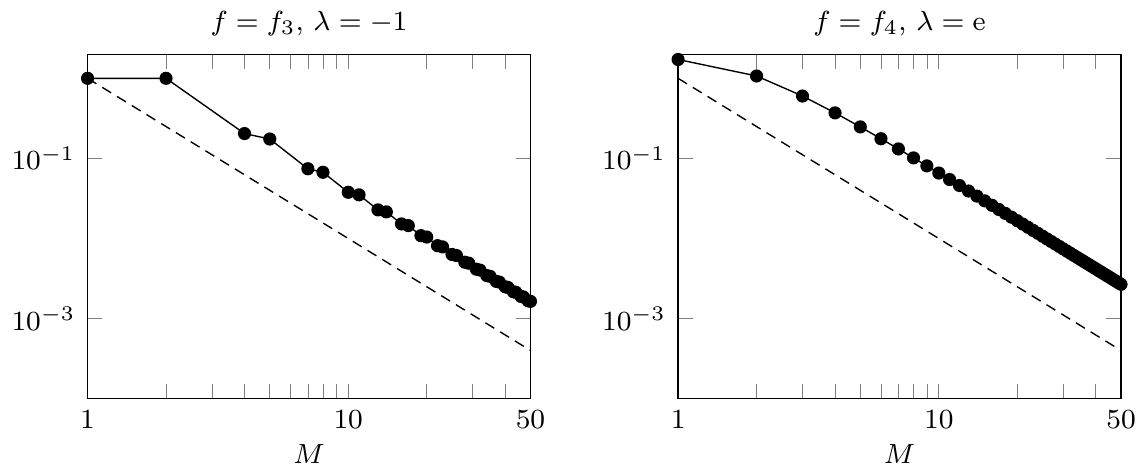}

\caption{Errors on $\lambda$ in the spectrum of $U_{\mathbb{C}}$ for $f = f_3$ in \cref{f3} and $f = f_4$ in \cref{f4}, computed with $L=1$.
The errors are the absolute errors on the approximated eigenvalue closest to $\lambda$.
Reference dashed lines show $M^{-2}$.}
\label{fig:spectra1err}
\end{figure}

In the examples in this \namecref{sec:experiments}, the isolated elements are actually determined exactly, at least with $M$ not too small, as in the cases of $f=f_1$ and $f=f_2$; some points very close to them, but not exactly them, are sometimes present as well.
The same happens for $1$ in the case of $f = f_4$.
The other elements on the boundaries, i.e., $e$ for $f=f_4$, as mentioned above, and $1$ and $-1$ for $f=f_3$, are approximated with order $2$, see \cref{fig:spectra1err}.
The finite order of convergence is a potential manifestation of non-smooth eigenfunctions; in fact, for smooth eigenfunctions the infinite order of pseudospectral methods is proved for DDEs \cite{BredaMasetVermiglio2012} and REs proper \cite{BredaLiessi2018}.

\begin{figure}
\centering
\includegraphics{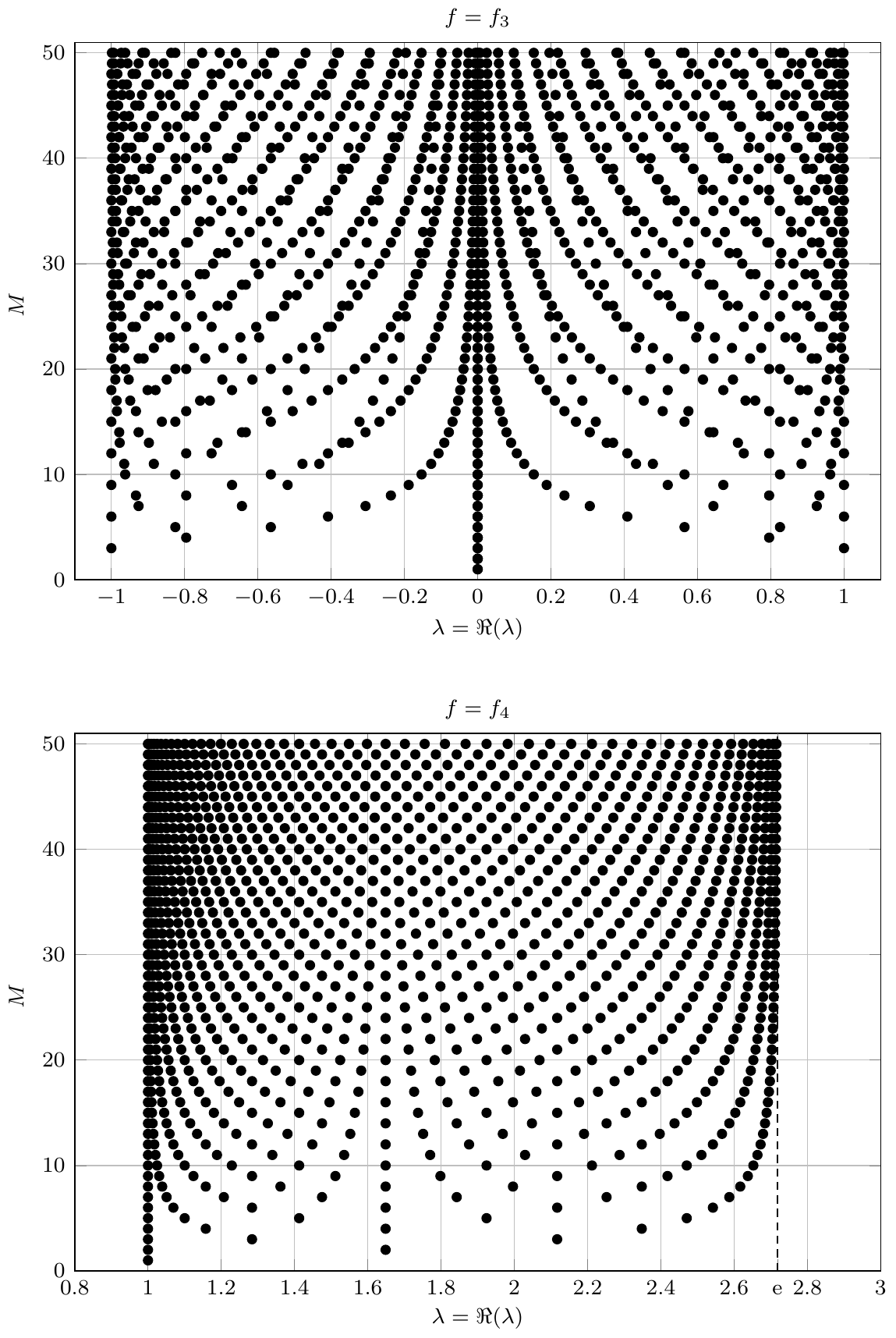}

\caption{Spectra of $U_{\mathbb{C}}$ with $f = f_3$ in \cref{f3} and $f = f_4$ in \cref{f4}, computed with $L=1$ and varying $M$.
Recall that the spectra are real.}
\label{fig:spectraM}
\end{figure}

As for the points in the interior, in order to have an intuition on the possible meaning of the convergence, in \cref{fig:spectraM} we provide plots of the spectra for $f = f_3$ and $f = f_4$: as $M$ increases, the overall distribution of the approximated spectrum seems to become denser.
We can observe that different points seem to appear in the approximated spectra periodically with respect to $M$.
This is true also for $1$ and $-1$ in the case of $f = f_3$, which may suggest that the order of convergence determined above does not describe the phenomenon accurately; for $e$ in the case of $f = f_4$, instead, it seems clear that the order of convergence makes sense.

\section{Linear systems with one delay}
\label{sec:system}

We consider now the linear NRE
\begin{equation}
\label{nre-system}
x(t) = A(t) x(t-\tau)
\end{equation}
with $A\colon \mathbb{R} \to \mathbb{R}^{d \times d}$ for $d\geq 2$ a $\tau$-periodic function with bounded variation components continuous from the right and we define $U$ as in \cref{sec:equation}.
We compute the spectrum of $U_{\mathbb{C}}$ for different choices of $A$ with $d=2$.

If $A$ is constant, the computed spectrum approximates $\sigma(A) \cup \{0\}$ even for small values of $M$, with some points of the spectrum being the exact ones and the others being very close to the exact ones.
Examples of matrices we used are
\begin{gather*}
\begin{pmatrix}
1 & 0 \\
0 & 1
\end{pmatrix}
,\quad
\begin{pmatrix}
1 & 1 \\
0 & 1
\end{pmatrix}
,\quad
\begin{pmatrix}
2 & 0 \\
0 & 1
\end{pmatrix}
,\quad
\begin{pmatrix}
2 & 1 \\
0 & 1
\end{pmatrix}
.
\end{gather*}

\begin{figure}
\centering
\includegraphics{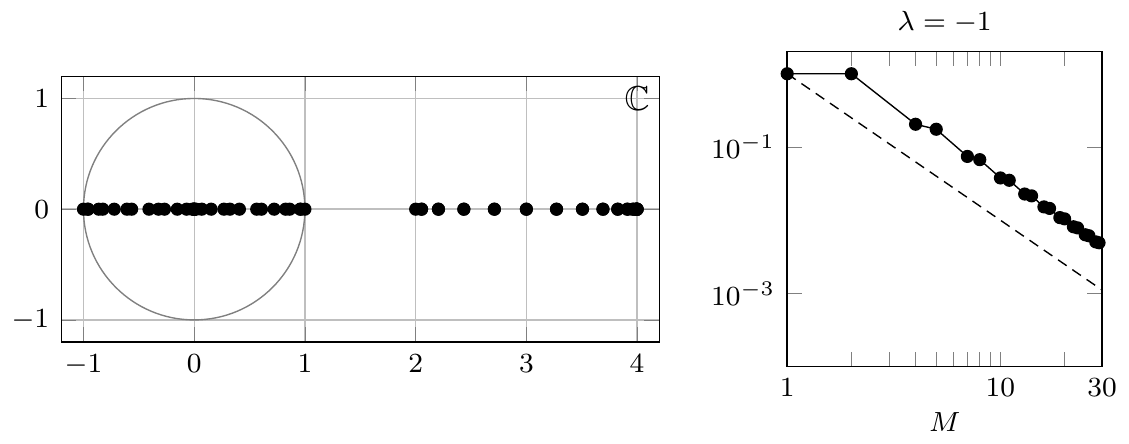}

\caption{Spectrum of $U_{\mathbb{C}}$ for \cref{nre-system} with $A$ defined by \cref{matrixstar} computed with $L=1$ and $M = 30$ (left) and errors on $\lambda = -1$ varying $M$ (right).
The errors are the absolute errors on the approximated eigenvalue closest to $\lambda$.
The reference dashed line shows $M^{-2}$.}
\label{fig:spectra-sys}
\end{figure}

For the experiments with a periodic non-constant $A$ we used matrices of the form
\begin{equation}
\label{matrixstar}
\begin{pmatrix}
\sin(2\pi t) & * \\
0 & 3+\cos(2\pi t)
\end{pmatrix}
\end{equation}
with constants or $1$-periodic functions in place of $*$.
The computed spectra approximate the union of the images of the diagonal elements, i.e., $[-1, 1] \cup [2, 4]$.
The elements $0$ and $4$ are computed exactly; the approximations of $-1$, $1$ and $2$ converge with order $2$.
See \cref{fig:spectra-sys} for an example with $\exp(t-\lfloor t\rfloor)-5$ in place of $*$.

\begin{figure}
\centering
\includegraphics{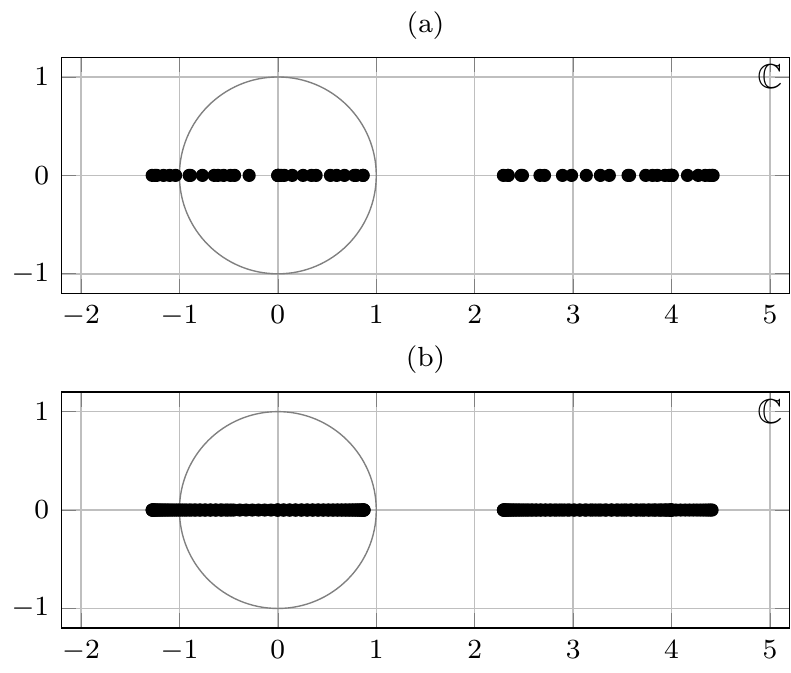}
\par
\includegraphics{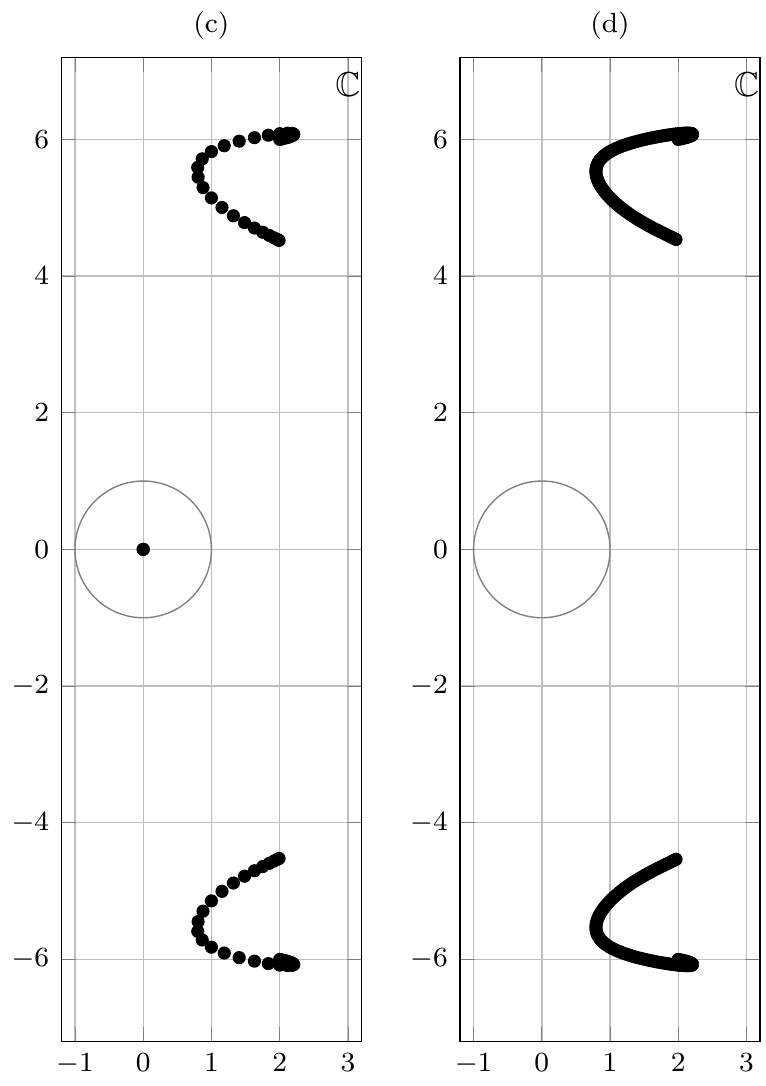}

\caption{Spectra of $U_{\mathbb{C}}$ computed with $L=1$ and $M=30$ (a,c), compared with the union of the spectra of $A(t)$ for $t$ varying on a uniform grid of $100$ points in $[0, 1]$ (b,d), for $A(t) = A_1(t)$ (a,b) and $A(t) = A_2(t)$ (c,d).}
\label{fig:spectra2}
\end{figure}

We also used the matrices
\begin{gather*}
A_1(t) \coloneqq
\begin{pmatrix}
\sin(2\pi t) & \exp(t-\lfloor t\rfloor) \\
\log(1+\abs{t-\lfloor t\rfloor}) & 3+\cos(2\pi t)
\end{pmatrix}
,\\
A_2(t) \coloneqq
\begin{pmatrix}
\sin(2\pi t) & \exp(t-\lfloor t\rfloor)-5 \\
\log(1+\abs{t-\lfloor t\rfloor})+10 & 3+\cos(2\pi t)
\end{pmatrix}
.
\end{gather*}
\Cref{fig:spectra2} compares the spectra of $U_{\mathbb{C}}$ with the union of the spectra of $A_i(t)$ varying $t\in[0,1]$ in these two cases.
The spectra of $A(t)$ are computed on a uniform grid of $100$ points in $[0, 1]$.

The experiments suggest the following conjecture.

\begin{conjecture}
Let $A\colon \mathbb{R} \to \mathbb{R}^{d \times d}$ be a $\tau$-periodic function with bounded variation components continuous from the right, and let $U$ be the monodromy operator of \cref{nre-system} (in the sense specified in \cref{sec:equation} in case $f$ is constant).
The spectrum of $U_{\mathbb{C}}$ is
\begin{equation*}
\sigma(U_{\mathbb{C}}) = \overline{\bigcup_{t \in \mathbb{R}} \sigma(A(t))}.
\end{equation*}
\end{conjecture}

\section{Linear scalar equations with two delays}
\label{sec:2delays}

We consider now the linear scalar NRE
\begin{equation}
\label{nre-2delays}
x(t) = a(t) x\Bigl(t-\frac{\tau}{2}\Bigr) + b(t) x(t-\tau)
\end{equation}
with $a, b \colon \mathbb{R} \to \mathbb{R}$ $\tau$-periodic functions of bounded variation continuous from the right and again we define $U$ as in \cref{sec:equation}.

\begin{theorem}
\label{th:2delays}
Let $\lambda = \alpha + \ii\beta \in \mathbb{C}$, let $a$ and $b$ be as described above and let $U$ be the monodromy operator of \cref{nre-2delays} (in the sense specified in \cref{sec:equation} in case $f$ is constant).
The operator $U_{\mathbb{C}} - \lambda I_{Y_{\mathbb{C}}}$ is bijective if and only if the determinant of
\begin{equation*}
\begin{pmatrix}
a(\theta) a(\theta+\frac{\tau}{2}) + b(\theta) -\alpha & \beta & a(\theta) b(\theta+\frac{\tau}{2}) & 0 \\
-\beta & a(\theta) a(\theta+\frac{\tau}{2}) + b(\theta) - \alpha & 0 & a(\theta) b(\theta+\frac{\tau}{2}) \\
a(\theta-\frac{\tau}{2}) & 0 & b(\theta-\frac{\tau}{2}) - \alpha & \beta \\
0 & a(\theta-\frac{\tau}{2}) & - \beta & b(\theta-\frac{\tau}{2}) - \alpha
\end{pmatrix}
\end{equation*}
is bounded away from $0$ for all $\theta \in [-\frac{\tau}{2}, 0]$.
In particular, if $a$ and $b$ are constant, $\lambda \in \mathbb{R}$ and $a^4+4a^2b\geq 0$, then
\begin{equation*}
\frac{a^2+2b\pm\sqrt{a^4+4a^2b}}{2} \in \sigma_{p}(U_{\mathbb{C}}).
\end{equation*}
\end{theorem}
\begin{proof}
Observe first that for $\theta \in [-\tau, 0)$ (recall that the elements of $Y$ have value $0$ at $0$), using the periodicity of $a$ and $b$,
\begin{equation*}
\begin{aligned}
(U\phi)(\theta) &= x_{\tau}(\theta; \phi) = x(\tau+\theta; \phi) = a(\tau+\theta) x\Bigl(\frac{\tau}{2}+\theta; \phi\Bigr) + b(\tau+\theta) x(\theta; \phi) \\
&=
\begin{cases}
a(\theta) (a(\frac{\tau}{2}+\theta) \phi(\theta) + b(\frac{\tau}{2}+\theta) \phi(\theta-\frac{\tau}{2})) + b(\theta) \phi(\theta) & \text{if $\theta \in [-\frac{\tau}{2}, 0)$,} \\
a(\theta) \phi(\frac{\tau}{2}+\theta) + b(\theta) \phi(\theta) & \text{if $\theta \in [-\tau, -\frac{\tau}{2})$.}
\end{cases}
\end{aligned}
\end{equation*}

Consider the equation $(U_{\mathbb{C}} - \lambda I_{Y_{\mathbb{C}}}) \phi = \psi$ for $\phi, \psi \in Y_{\mathbb{C}}$, equivalent to $U\phi_{\Re} + \ii U\phi_{\Im} - \lambda\phi_{\Re} - \ii\lambda\phi_{\Im} = \psi_{\Re} + \ii\psi_{\Im}$ for $(\phi_{\Re}, \phi_{\Im}), (\psi_{\Re}, \psi_{\Im}) \in Y \times Y$ by separating the real and imaginary parts.
The latter is equivalent to
\begin{equation*}
\begin{cases}
a(\theta) (a(\frac{\tau}{2}+\theta) \phi_{\Re}(\theta) + b(\frac{\tau}{2}+\theta) \phi_{\Re}(\theta-\frac{\tau}{2})) + b(\theta) \phi_{\Re}(\theta) \\
+ \ii(a(\theta) (a(\frac{\tau}{2}+\theta) \phi_{\Im}(\theta) + b(\frac{\tau}{2}+\theta) \phi_{\Im}(\theta-\frac{\tau}{2})) + b(\theta) \phi_{\Im}(\theta)) \\
- (\alpha+\ii\beta)\phi_{\Re}(\theta) - \ii(\alpha+\ii\beta)\phi_{\Im}(\theta) = \psi_{\Re}(\theta) + \ii\psi_{\Im}(\theta) & \text{if $\theta \in [-\frac{\tau}{2}, 0)$,} \\[1ex]
a(\theta) \phi_{\Re}(\frac{\tau}{2}+\theta) + b(\theta) \phi_{\Re}(\theta)
+ \ii(a(\theta) \phi_{\Im}(\frac{\tau}{2}+\theta) + b(\theta) \phi_{\Im}(\theta)) \\
- (\alpha+\ii\beta)\phi_{\Re}(\theta) - \ii(\alpha+\ii\beta)\phi_{\Im}(\theta) = \psi_{\Re}(\theta) + \ii\psi_{\Im}(\theta) & \text{if $\theta \in [-\tau, -\frac{\tau}{2})$,}
\end{cases}
\end{equation*}
and, by separating the real and imaginary parts, to
\begin{equation*}
\begin{cases}
a(\theta) a(\frac{\tau}{2}+\theta) \phi_{\Re}(\theta) + a(\theta) b(\frac{\tau}{2}+\theta) \phi_{\Re}(\theta-\frac{\tau}{2}) + b(\theta) \phi_{\Re}(\theta) \\
- \alpha \phi_{\Re}(\theta) + \beta \phi_{\Im}(\theta) = \psi_{\Re}(\theta) & \text{if $\theta \in [-\frac{\tau}{2}, 0)$,} \\[1ex]
a(\theta) a(\frac{\tau}{2}+\theta) \phi_{\Im}(\theta) + a(\theta) b(\frac{\tau}{2}+\theta) \phi_{\Im}(\theta-\frac{\tau}{2}) + b(\theta) \phi_{\Im}(\theta) \\
- \beta \phi_{\Re}(\theta) - \alpha \phi_{\Im}(\theta) = \psi_{\Im}(\theta) & \text{if $\theta \in [-\frac{\tau}{2}, 0)$,} \\[1ex]
a(\theta) \phi_{\Re}(\frac{\tau}{2}+\theta) + b(\theta) \phi_{\Re}(\theta) - \alpha \phi_{\Re}(\theta) + \beta \phi_{\Im}(\theta) = \psi_{\Re}(\theta) & \text{if $\theta \in [-\tau, -\frac{\tau}{2})$,} \\[1ex]
a(\theta) \phi_{\Im}(\frac{\tau}{2}+\theta) + b(\theta) \phi_{\Im}(\theta) - \beta \phi_{\Re}(\theta) - \alpha \phi_{\Im}(\theta) = \psi_{\Im}(\theta) & \text{if $\theta \in [-\tau, -\frac{\tau}{2})$,}
\end{cases}
\end{equation*}
and in turn, with a change of variable in the two latter equations, to
\begin{equation*}
\begin{cases}
a(\theta) a(\frac{\tau}{2}+\theta) \phi_{\Re}(\theta) + a(\theta) b(\frac{\tau}{2}+\theta) \phi_{\Re}(\theta-\frac{\tau}{2}) + b(\theta) \phi_{\Re}(\theta) \\
- \alpha \phi_{\Re}(\theta) + \beta \phi_{\Im}(\theta) = \psi_{\Re}(\theta) & \text{if $\theta \in [-\frac{\tau}{2}, 0)$,} \\[1ex]
a(\theta) a(\frac{\tau}{2}+\theta) \phi_{\Im}(\theta) + a(\theta) b(\frac{\tau}{2}+\theta) \phi_{\Im}(\theta-\frac{\tau}{2}) + b(\theta) \phi_{\Im}(\theta) \\
- \beta \phi_{\Re}(\theta) - \alpha \phi_{\Im}(\theta) = \psi_{\Im}(\theta) & \text{if $\theta \in [-\frac{\tau}{2}, 0)$,} \\[1ex]
a(\theta-\frac{\tau}{2}) \phi_{\Re}(\theta) + b(\theta-\frac{\tau}{2}) \phi_{\Re}(\theta-\frac{\tau}{2}) \\
- \alpha \phi_{\Re}(\theta-\frac{\tau}{2}) + \beta \phi_{\Im}(\theta-\frac{\tau}{2}) = \psi_{\Re}(\theta-\frac{\tau}{2}) & \text{if $\theta \in [-\frac{\tau}{2}, 0)$,} \\[1ex]
a(\theta-\frac{\tau}{2}) \phi_{\Im}(\theta) + b(\theta-\frac{\tau}{2}) \phi_{\Im}(\theta-\frac{\tau}{2}) \\
- \beta \phi_{\Re}(\theta-\frac{\tau}{2}) - \alpha \phi_{\Im}(\theta-\frac{\tau}{2}) = \psi_{\Im}(\theta-\frac{\tau}{2}) & \text{if $\theta \in [-\frac{\tau}{2}, 0)$.}
\end{cases}
\end{equation*}
We can write the latter in matrix form for $\theta \in [-\frac{\tau}{2}, 0)$ as
\begin{equation*}
\begin{aligned}
\begin{pmatrix}
a(\theta) a(\theta+\frac{\tau}{2}) + b(\theta) -\alpha & \beta & a(\theta) b(\theta+\frac{\tau}{2}) & 0 \\
-\beta & a(\theta) a(\theta+\frac{\tau}{2}) + b(\theta) - \alpha & 0 & a(\theta) b(\theta+\frac{\tau}{2}) \\
a(\theta-\frac{\tau}{2}) & 0 & b(\theta-\frac{\tau}{2}) - \alpha & \beta \\
0 & a(\theta-\frac{\tau}{2}) & - \beta & b(\theta-\frac{\tau}{2}) - \alpha
\end{pmatrix}
\\
\cdot
\begin{pmatrix}
\phi_{\Re}(\theta) \\
\phi_{\Im}(\theta) \\
\phi_{\Re}(\theta-\frac{\tau}{2}) \\
\phi_{\Im}(\theta-\frac{\tau}{2})
\end{pmatrix}
=
\begin{pmatrix}
\psi_{\Re}(\theta) \\
\psi_{\Im}(\theta) \\
\psi_{\Re}(\theta-\frac{\tau}{2}) \\
\psi_{\Im}(\theta-\frac{\tau}{2})
\end{pmatrix}
\end{aligned}
\end{equation*}
and the first part of the thesis follows by observing that according to the Laplace formula the inverse of the matrix is the adjugate matrix divided by the determinant.

As for the second part, if $a$ and $b$ are constant and $\lambda \in \mathbb{R}$, i.e., $\beta = 0$, the matrix becomes
\begin{equation*}
\begin{pmatrix}
a^2 + b -\alpha & 0 & a b & 0 \\
0 & a^2 + b - \alpha & 0 & a b \\
a & 0 & b - \alpha & 0 \\
0 & a & 0 & b - \alpha
\end{pmatrix}
.
\end{equation*}
Its determinant is $((a^2+b-\alpha)(b-\alpha)-a^2b)^2$, which is $0$ (and thus $U_{\mathbb{C}} - \lambda I_{Y_{\mathbb{C}}}$ is not injective) if and only if $\alpha^2-(a^2+2b)\alpha+b^2=0$, i.e., if and only if $a^4+4a^2b\geq0$ and
\begin{equation*}
\lambda = \alpha = \frac{a^2+2b\pm\sqrt{a^4+4a^2b}}{2}.
\qedhere
\end{equation*}
\end{proof}

\begin{figure}
\centering
\includegraphics{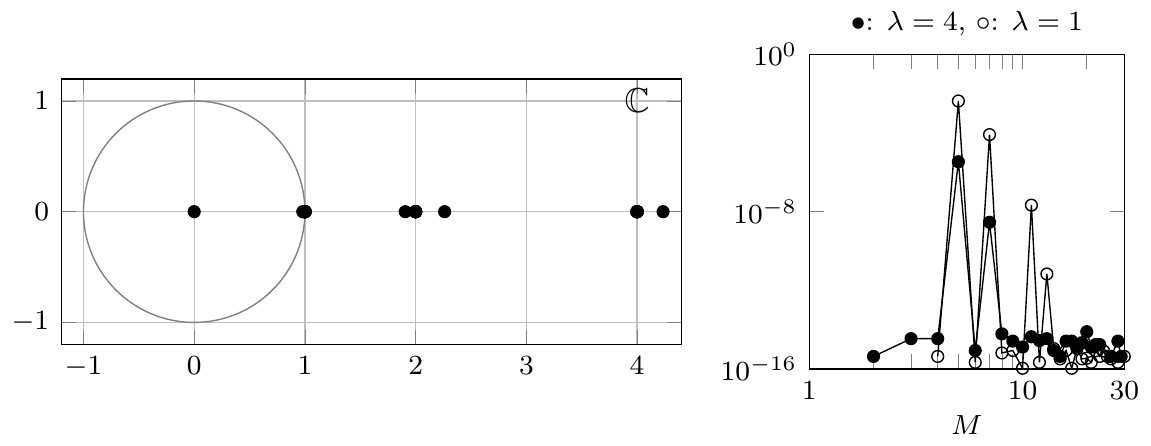}

\caption{Spectrum of $U_{\mathbb{C}}$ for \cref{nre-2delays} with $a \equiv 1$ and $b \equiv 2$ computed with $L=1$ and $M = 30$ (left) and errors on the known points $\lambda \in \sigma(U_{\mathbb{C}})$ (as shown in the legend) varying $M$ (right).
The errors are the absolute errors on the approximated eigenvalue closest to $\lambda$.
Missing dots correspond to null error.}
\label{fig:spectra-2delays-err}
\end{figure}

\Cref{fig:spectra-2delays-err} shows the spectrum of $U_{\mathbb{C}}$ computed with $L=1$ and $M=30$ and the errors on its known real elements given by \cref{th:2delays} for $a\equiv 1$ and $b\equiv2$, i.e., $1$ and $4$: depending on $M$ the approximation error is of the order of the machine precision or is apparently decaying with infinite order.
The approximated spectrum (also varying $M$) seems to suggest that in this case $\sigma(U_{\mathbb{C}})$ contains only isolated points.

\section{Finally, a perspective on convergence}
\label{sec:conclusions}

The aim of this work has been to begin a theoretical and numerical investigation on the spectra of evolution operators of NREs, with the stability of equilibria and periodic orbits in mind.
We already noted that the theory linking the latter to the former is currently lacking, although it seems plausible in light of \cite{DiekmannVerduynLunel2021}.
We also noted that we did not attempt to prove the well-posedness of the discretized collocation equation \cref{fixed-point-UMN} and the convergence of the numerical method (not to mention that the precise meaning of convergence in this context still needs to be specified).
Indeed, the proofs of those results for REs in \cite{BredaLiessi2018} are based on the Banach Perturbation Lemma and require that the equation has a regularizing effect, which is not the case for NREs.

A possible alternative approach to the convergence of the method is to consider REs or DDEs that in some sense approximate the NREs, and to apply to them the methods of \cite{BredaLiessi2018,BredaMasetVermiglio2012}, for which the convergence has been proved.
Comparing the results may allow us both to assess the behavior of the numerical method for NREs and to understand the relation between these equations.

With reference to \cref{nre}, given $\epsilon>0$, we may consider the linear RE
\begin{equation}
\label{relim}
x(t) = \frac{f(t)}{\epsilon} \int_{-\tau}^{-\tau+\epsilon} x(t + \theta) \dd \theta
\end{equation}
and the linear DDE
\begin{equation}
\label{ddelim}
y'(t) = \frac{f(t)}{\epsilon} (y(t - \tau + \epsilon) - y(t - \tau)),
\end{equation}
Observe that solutions $y$ of \cref{ddelim} are primitives of solutions $x$ of \cref{relim}.

Let $U$ be the monodromy operator of \cref{nre}, as above, and let $U_{\epsilon}$ be that of \cref{relim} (similar arguments apply to \cref{ddelim}).
It can be proved that there is a pointwise convergence of the solutions of \cref{relim} to those of \cref{nre}, i.e., for each $\phi\in Y$ and $\theta \in [-\tau, 0]$ we can prove that
\begin{equation}
\label{Uepsconv}
\lim_{\epsilon \to 0^+} (U_{\epsilon}\phi)(\theta) = (U\phi)(\theta).
\end{equation}
However, it is also true that $U_{\epsilon}$ is compact (or has a compact power) \cite{BredaLiessi2021}, while neither $U$ nor its powers are, which casts doubts on how strong the meaning of the convergence of \cref{relim} to \cref{nre} (and of the corresponding monodromy operator and spectra) can be.

Let $U_{M}$ and $U_{\epsilon,M}$ be the discretized versions of $U$ and $U_{\epsilon}$ with $L=1$.
We know from \cite{BredaLiessi2018} that as $M$ increases the spectrum of $U_{\epsilon,M}$ converges to that of $U_{\epsilon}$ in the sense specified therein.
We also showed in the experiments of \cref{sec:experiments} that the spectrum of $U_{M}$ converges to that of $U$ in the intuitive sense specified above.
Moreover, in other experiments not presented here we see that the spectrum of $U_{\epsilon,M}$ converges to that of $U_{M}$ as $\epsilon$ vanishes.
The experimental evidence seems thus to validate the idea that the convergence of $U_{\epsilon}$ to $U$ is stronger than \cref{Uepsconv} and that the spectrum of the former should converge in some sense to that of the latter.

The investigation of this topic will be the subject of future work, which may be informed by the results in \cite{ChowDiekmannMalletParet1985}, where the properties of an equation similar to \cref{relim} are studied, and in \cite{MalletParetNussbaum1986,MalletParetNussbaum1989}, which deal with an equation with some similarities to \cref{ddelim}.

\section*{Acknowledgments}

Dimitri Breda and Davide Liessi are members of INdAM Research group GNCS and of UMI Research group ``Mo\-del\-li\-sti\-ca socio-epidemiologica''.

This work was partially supported by the Italian Ministry of University and Research (MUR) through the PRIN 2020 project (No.\ 2020JLWP23) ``Integrated Mathematical Approaches to Socio-Epidemiological Dynamics'' (CUP: E15F21005420006).
The work of Davide Liessi was partially supported by Finanziamento Giovani Ricercatori 2018--2019 and 2020--2021 of INdAM Research group GNCS.

{\sloppy
\printbibliography
\par}

\end{document}